\documentclass[12pt,reqno]{amsart}
\usepackage{amsfonts,amsmath,amssymb}
\usepackage{tikz}
\usetikzlibrary{decorations.pathmorphing}

\usepackage{graphicx}

\def\bim{\begin{itemize}\item[]}
	
	\def\eim{\end{itemize}}

\def\uhat{\widehat{u}}
\def\ubar{\overline{u}}

\allowdisplaybreaks
\advance\textheight by \topskip

\usepackage{mathrsfs}

\usepackage[latin1]{inputenc}
\usepackage{graphicx}

\newtheorem{theorem}{Theorem}
\usepackage{euscript}

\def\[{[\! [}
\def\]{]\! ]}

\title[MIN-turns and MAX-turns in $k$-Dyck paths]{MIN-turns and MAX-turns in $k$-Dyck paths: a pure generating function approach}

\author[H.~Prodinger]{Helmut Prodinger}

	\address{Department of Mathematics, University of Stellenbosch 7602, Stellenbosch, South Africa
	and
	NITheCS (National Institute for
	Theoretical and Computational Sciences), South Africa.}
\email{hproding@sun.ac.za}

\date{\today}

\begin{document}

\begin{abstract}
	$k$-Dyck paths differ from ordinary Dyck paths by using an up-step of length $k$. We analyze at which level the path is after the $s$-th up-step and
	before the $(s+1)$st up-step. In honour of Rainer Kemp who studied a related concept 40 years ago the terms \textsc{max}-terms and \textsc{min}-terms
	are used. Results are obtained by an appropriate use of trivariate generating functions; practically no combinatorial arguments are used.
\end{abstract}

\maketitle

\section{Introduction}

The objects of our interest in this paper are $k$-Dyck paths, having up-steps $(1,k)$ and down-steps $(1,-1)$, and never go below the $x$-axis.
At the end, they reach the $x$-axis, but we also need versions that end at a prescribed level different from 0.
Much material about such paths can be found in \cite{selkirk-master}, \cite{garden} and in many standard textbooks on combinatorics.

We consider \textsc{max}-turns, where each up-step ends and \textsc{min}-turns, where each up-step starts.
The Figure~\ref{basic-fig} explains the concept readily:
\begin{figure}[h]\label{basic-fig}
\begin{center}
	\begin{tikzpicture}[scale=0.5]
		\draw (0,0) -- (30,0);
		\draw (0,0) -- (0,13);
		\draw[thick](0,0)--(1,5);
		\draw [decorate,decoration=snake,thick] (1,5) -- (3,3);
		\draw[thick](3,3)--(4,8);
		\draw [decorate,decoration=snake,thick] (4,8) -- (5,7);
		\draw[thick](5,7)--(6,12);
		\draw [decorate,decoration=snake,thick] (6,12) -- (10,8);
		\draw[thick](10,8)--(11,13);
		\draw [decorate,decoration=snake,thick] (11,13) -- (18,6);
		
				\draw [decorate,decoration=snake,thick,brown] (18,6) to [in=100,out=60] (30,0);
		
		\node[thick, cyan] at (3,3){$\bullet$};
		\node[thick, cyan] at (5,7){$\bullet$};
		\node[thick, red] at (4,8){$\bullet$};
		\node[thick, red] at (6,12){$\bullet$};
		\node[thick, red] at (1,5){$\bullet$};
		\node[thick, red] at (11,13){$\bullet$};
		\node[thick, cyan] at (10,8){$\bullet$};
				\node[thick, cyan] at (18,6){$\bullet$};
		

	\end{tikzpicture}
\end{center}
\caption{The first 4 \textsc{max}-turns and the first 4 \textsc{min}-turns are shown.}
\end{figure}
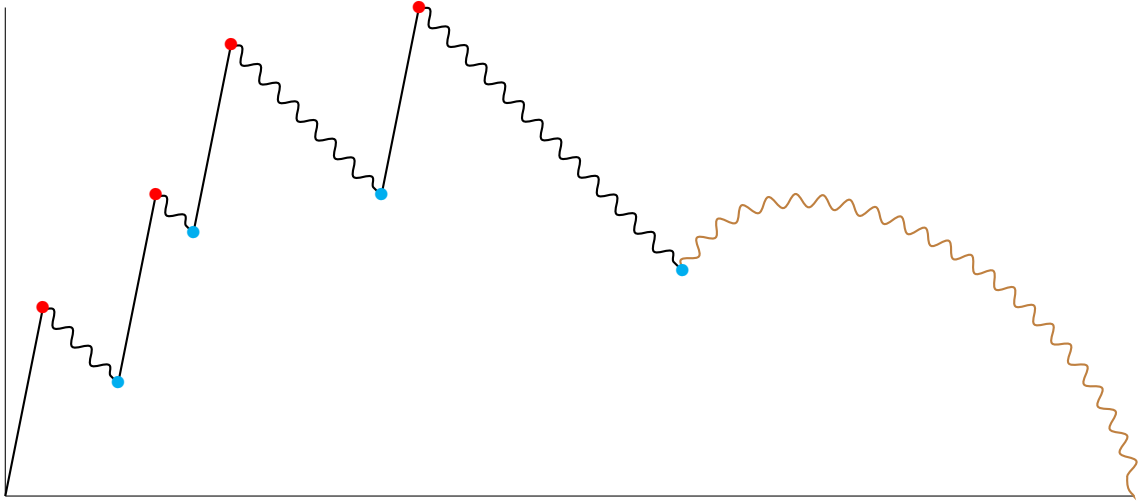	

The $k$-Dyck paths can only exist for a length of the form $(k+1)N$, which is clear for combinatorial reasons or otherwise.
We want to know the average level of the $s$-th \textsc{max}-turn resp.\ \textsc{min}-turn, amoung all $k$-Dyck of the same length.
In order to do this, we sum the level of the $s$-th \textsc{max}-turn resp.\ \textsc{min}-turn, over all $k$-Dyck paths of  length $(k+1)N$.
The generating functions we obtain, enumerate the (partial) $k$-Dyck paths that end at a prescribed level $j$ as well; this is of interest in \cite{deng} and
later in \cite{garden}, in the section on Hoppy's walks.  

To get the average, one only needs to divide by the total number of such $k$-Dyck paths.
Our main achievement is to get the \emph{fully explicit} functions
\begin{align*}
\textsc{max}(z,w)=\sum_{N\ge0,\, s\ge1}z^{(k+1)N}w^s[&\text{cumulative level of the $s$-th \textsc{max}-turn}\\[-15pt]
& \text{in all $k$-Dyck paths of length $(k+1)N$}]
\end{align*}
and
\begin{align*}
	\textsc{min}(z,w)=\sum_{N\ge0,\, s\ge1}z^{(k+1)N}w^s[&\text{cumulative level of the $s$-th \textsc{min}-turn}\\[-15pt]
	& \text{in all $k$-Dyck paths of length $(k+1)N$}].
\end{align*}
As a bonus we get $\textsc{osc}(z,w):=\textsc{max}(z,w)-\textsc{min}(z,w)$; this cumulates the lengths of the wavy line between
the $s$-th \textsc{max}-turn and $s$-th \textsc{min}-turn, see Figure~\ref{basic-fig}. The function $\textsc{osc}(z,w)$  is  somewhat simpler than 
$\textsc{min}(z,w)$ and $\textsc{max}(z,w)$, because of cancellations. In this way, we recover some of the results
from \cite{AHS} without resorting to any bijective combinatorics. Note that such a wavy line might have length zero as well, if two
up-steps follow each other immediately. Rainer Kemp in \cite{Kemp} has considered \textsc{max}- and \textsc{min}-turns for Dyck paths,
although his definitions were slightly different (peaks and valleys).

The key to the success of our method is the simple but perhaps unusual identity
\begin{equation*}
\sum_{i\ge0}\big([z^i]f(z)\big)\cdot y^i=f(y).
\end{equation*}

\section{Some basic observations}

As we will see soon, the equation $u=z+zwu^{k+1}$ plays a major role when enumerating $k$-Dyck paths. The equation is of the form
$u=z\Phi(u)$, with $\Phi(u)=1+wu^{k+1}$, so it is amenable to the Lagrange inversion \cite{FS}, and the coefficients of the inverse function 
$u=u(z)$ can be computed:
\begin{align*}
[z^{(k+1)N+1}]u&=\frac1{(k+1)N+1}[u^{(k+1)N}](1+wu^{k+1})^{(k+1)N+1}\\&=\frac1{(k+1)N+1}w^N\binom{(k+1)N+1}{N}.
\end{align*}
One solution of the equation, which is analytic at the origin and has combinatorial significance, is called $\ubar$. We will use the formula 
(again by the Lagrange inversion formula)
\begin{equation}\label{ubarnumber}
\ubar=z\sum_{N\ge0}w^Nz^{(k+1)N}\frac1{kN+1}\binom{(k+1)N}{N}.
\end{equation}
We also need the specialization where $w=1$, which we call $\uhat$, which satisfies the equation $u=z+zu^{k+1}$:
\begin{equation}\label{uhatnumber}
	\uhat=z\sum_{\lambda\ge0}\frac{1}{(k+1)\lambda+1}\binom{1+(k+1)\lambda}{\lambda}z^{(k+1)\lambda}
\end{equation}
and (again by the Lagrange inversion formula)
\begin{equation}\label{minusknumber}
	\uhat^{-k}=z^{-k}-z\sum_{\lambda\ge0}\frac{k}{\lambda+1}\binom{(k+1)\lambda}{\lambda}z^{(k+1)\lambda}.
\end{equation}

\section{\textsc{min}-turns}

After these classical observations about the implicit function, we move to the treatment of \textsc{min}-turns.
Our method of choice is \emph{adding a new slice}.
The following substitution is essential for adding a new slice (which is an up-step, followed by a maximal sequence of down-steps):
\begin{equation*}
	u^j\longrightarrow zw\sum_{0\le i \le j+k} z^iu^{j+k-i}=\frac{zwu^{k+1}}{u-z}u^j-\frac{wz^{k+2}}{u-z}z^j.
\end{equation*}
The substitution tells us what can become of a partial $k$-ary path landing at level $j$ when adding (attaching) another slice landing on level $j+k-i$.
The technique of adding-a-new slice is described in \cite{FP-slice}.

Now let $F_m(u)=F_m(u;z)$ be the generating function according to $m$ slices; $z$ refers to the lengths and $u$ to the level of the $m$-th \textsc{min}-turn. The substitution leads to
\begin{equation*}
	F_{m+1}(u)=\frac{zwu^{k+1}}{u-z}F_m(u)-\frac{wz^{k+2}}{u-z}F_m(z),\quad F_0(u)=1.
\end{equation*}
Let $F(u)=\sum_{m\ge0}F_m(u)$, so that we do not care about the number $m$ anymore, since the variable $w$ takes care of it; then
\begin{equation*}
	F(u)=F(u;z,w)=1+\frac{zwu^{k+1}}{u-z}F(u)-\frac{wz^{k+2}}{u-z}F(z),
\end{equation*}
or
\begin{equation*}
	F(u)=\frac{u-z-wz^{k+2}F(z)}{u-z-zwu^{k+1}}.
\end{equation*}	
$u=\ubar$ is a factor of the denominator, but for $z$ and $u$ small, we have $u\sim z$, so this factor must cancel in the numerator as well. This is what
one learns from the kernel method \cite{prodinger-kernel}.
We find
\begin{equation*}
F(z)=\frac{\overline{u}-z}{wz^{k+2}}
\end{equation*}
and further
\begin{equation*}
	F(u)=F(u;z,w)=\frac{u-\overline{u}}{u-z-zwu^{k+1}}.
\end{equation*}	
That ends the computation of the ``left'' part of the $k$-Dyck path. For the right one,  we start at level $h$ with an up-step and end eventually at the zero level.
We don't use the variable $w$ here. The kernel method could also be used, but there is a simpler way. 
Reading the path from right to left, there is a decomposition when the path leaves a level and never comes back to it.
Recall  that $\ubar/z$ is the generating function of $k$-Dyck paths.
With this, we find for the enumeration of the ``right'' part
\begin{equation*}
\Big(\frac\ubar z\Big)^hz^h\Big(\frac\ubar z-1\Big)=\frac{(\overline{u}-z)\overline{u}^h}{z};
\end{equation*}
the minus 1 term happens since the reversed path must end with a down-step. The formula works only for $h\ge1$, but the
instance $h=0$ is not needed (although easy).
The fact that this generating function is essentially a \emph{power} is part of our successful approach.
We now have
\begin{align*}
\textsc{min}(z,w)&=\sum_{h\ge1}h[u^h]F(u)\cdot \frac{(\overline{u}-z)\overline{u}^h}{z}\\
&=\sum_{h\ge1}[u^{h-1}]\frac{d}{du}F(u)\cdot \frac{(\overline{u}-z)\overline{u}^h}{z}\\
&=\frac{(\overline{u}-z)\overline{u}}{z}\cdot \frac{-zu+uzw{u}^{k+1}k+\ubar u-\ubar zw{u}^{k+1}k-\ubar zw{u}^{k+1}}{u(u-z-zwu^{k+1})^2}\bigg|_{u=\uhat}\\
&=\frac{(\overline{u}-z)\overline{u}}{z}\cdot \frac{-z\uhat+uzw{\uhat}^{k+1}k+\ubar \uhat-\ubar zw{\uhat}^{k+1}k-\ubar zw{\uhat}^{k+1}}{\uhat(\uhat-z-zw\uhat^{k+1})^2}.
\end{align*}
A simplification that only uses $\uhat^{k+1}=\frac{\uhat-z}{z}$ eventually leads to
\begin{align*}
	\textsc{min}(z,w)=	\frac{kw\uhat}
	{z(1-w)^2}+
	\frac{(\ubar-z)}
	{z^2(1-w)^2\uhat^{k}}-
	\frac{(k+1)\overline{u}w}
	{z(1-w)^2}
	\end{align*}
which can be checked easily, even by hand.
\begin{theorem}
	The generating funtion $\textsc{min}(z,w)$ where the coefficient of $z^{(k+1)N}w^s$ refers to the cumulative levels of the $s$-th \textsc{min}-turn, is given by
	\begin{align*}
		\textsc{min}(z,w)=	\frac{kw\uhat}
	{z(1-w)^2}+
	\frac{(\ubar-z)}
	{z^2(1-w)^2\uhat^{k}}-
	\frac{(k+1)\overline{u}w}
	{z(1-w)^2}.
\end{align*}
\end{theorem}
The next step is to expand this function:
	\begin{align*}
	[w^s]\textsc{min}(z,w)&=	[w^s]\frac{kw\uhat}
	{z(1-w)^2}-[w^s]	\frac{(k+1)\overline{u}w}	{z(1-w)^2}+[w^s]	\frac{(\ubar-z)}
	{z^2(1-w)^2\uhat^{k}}\\
	&=	\frac{sk\uhat}	{z}
	-(k+1)\sum_{i=0}^{s-1}	(s-i)[w^i]\frac{\overline{u}}{z}\\
	&\quad+\sum_{i=0}^{s}(s+1-i)[w^i]	\frac{z\sum_{N\ge0}w^Nz^{(k+1)N}\frac1{kN+1}\binom{(k+1)N}{N}-z}
	{z^2\uhat^{k}}\\
	&=	\frac{sk\uhat}	{z}
-(k+1)\sum_{i=0}^{s-1}	(s-i)[w^i]\frac{\overline{u}}{z}\\
&\quad+\sum_{i=0}^{s}(s+1-i)[w^i]	\frac{\sum_{N\ge1}w^Nz^{(k+1)N}\frac1{kN+1}\binom{(k+1)N}{N}}
{z\uhat^{k}}\\
	&=	\frac{sk\uhat}	{z}
	-(k+1)\sum_{i=0}^{s-1}	(s-i)z^{(k+1)i}\frac1{ki+1}\binom{(k+1)i}{i}\\
	&\quad+\sum_{i=1}^{s}(s+1-i)	\frac{1}
	{z\uhat^{k}}z^{(k+1)i}\frac1{ki+1}\binom{(k+1)i}{i}\\
	&=	\frac{sk\uhat}	{z}
	-(k+1)\sum_{i=0}^{s-1}	(s-i)z^{(k+1)i}\frac1{ki+1}\binom{(k+1)i}{i}\\
	&\quad+\sum_{i=1}^{s}(s+1-i)z^{(k+1)(i-1)}\frac1{ki+1}\binom{(k+1)i}{i}\\
	&\quad-\sum_{i=1}^{s}(s+1-i)	
	z^{(k+1)i}\frac1{ki+1}\binom{(k+1)i}{i}\sum_{\lambda\ge0}\frac{k}{\lambda+1}\binom{(k+1)\lambda}{\lambda}z^{(k+1)\lambda}.
\end{align*}
And now we read off the coefficient of $z^{(k+1)N}$; we assume that $N\ge s$, otherwise a path would not have an $s$-th \textsc{min}-turn:
	\begin{align*}
	[w^sz^{(k+1)N}]\textsc{min}(z,w)
	&=	sk \frac1{kN+1}\binom{(k+1)N}{N}\\
	&-\sum_{i=1}^{s}(s+1-i)	\frac1{ki+1}\binom{(k+1)i}{i}\frac{k}{(N-i)+1}\binom{(k+1)(N-i)}{(N-i)}.
\end{align*}
\begin{theorem}
	The sum of levels of  the $s$-th \textsc{min}-turns in all the $k$-Dyck paths of length $(k+1)N$ is given by
	\begin{align*}
			 \frac{sk}{kN+1}\binom{(k+1)N}{N}-\sum_{i=1}^{s}(s+1-i)	\frac1{ki+1}\binom{(k+1)i}{i}\frac{k}{(N-i)+1}\binom{(k+1)(N-i)}{(N-i)}.
	\end{align*}
\end{theorem}
This formula is fully explicit and covers all instances of $N\ge s$.

\section{\textsc{max}-turns}

It is easy to go from a  \textsc{min}-turn to the next \textsc{max}-turn, just by doing \emph{one} up-step. On the level of generating functions, this means
\begin{equation*}
G(u;z,w)=F(u;z,w)wzu^{k}.
\end{equation*}
The right part of a part (everything from a given \textsc{max}-turn until the end) is even easier than before, since level $h$ must be reached without further restriction. The result for the ``right'' side of the path is then
\begin{equation*}
\Big(\frac{\uhat}{z}\Big)^{h+1}z^h=\uhat^{h-1}\frac{\uhat^2}{z}.
\end{equation*}
After these preliminaries, we can engage into the big computation:
\begin{align*}
\textsc{max}(z,w)&=\sum_{h\ge1}h[u^h]G(u)\cdot \uhat^{h-1}\frac{\uhat^2}{z}
=\frac{\uhat^2}{z}\sum_{h\ge1}[u^{h-1}]\frac{d}{du}G(u)\cdot \uhat^{h-1}\\
&=	\frac{wk \uhat}{z(1-w)^2}
	-	\frac{wk\ubar }{z(1-w)^2}
	+\frac{w (\ubar-z)}{z^2\uhat^{k}(1-w)^2}
		-\frac{w^2 \ubar }{z(1-w)^2}.
\end{align*}
The same type of simplifications as before has been applied.

And now we go to the coefficients of this, in a similar style as it was done for the \textsc{min}-turns:
\begin{align*}
[w^s]	\textsc{max}(z,w)
	&=	s\frac{k \uhat}{z}
	-	\sum_{i=0}^{s-1}(s-i)[w^i]\frac{k\ubar }{z}\\&\quad
	+\frac1{z\uhat^{k}}\sum_{i=0}^{s-1}(s-i)[w^i]\frac{ (\ubar-z)}{z}
	-\sum_{i=0}^{s-2}(s-1-i)[w^i]\frac{\ubar }{z}\\
	&=	sk\sum_{N\ge0}z^{(k+1)N}\frac1{kN+1}\binom{(k+1)N}{N}\\&\quad
	-	\sum_{i=0}^{s-1}(s-i)z^{(k+1)i}\frac1{ki+1}\binom{(k+1)i}{i}\\&\quad
	+\sum_{i=1}^{s-1}(s-i)z^{(k+1)(i-1)}\frac1{ki+1}\binom{(k+1)i}{i}\\&\quad
	-\sum_{\lambda\ge0}\frac{k}{\lambda+1}\binom{(k+1)\lambda}{\lambda}z^{(k+1)\lambda}\sum_{i=1}^{s-1}(s-i)z^{(k+1)i}\frac1{ki+1}\binom{(k+1)i}{i}\\&	\quad
	-\sum_{i=0}^{s-2}(s-1-i)z^{(k+1)i}\frac1{ki+1}\binom{(k+1)i}{i}.
\end{align*}
And, again for $N\ge s$, we read off the coefficient of $z^{(k+1)N}$:
\begin{align*}
	[w^sz^{(k+1)N}]	\textsc{max}(z,w)
	&=	sk\frac1{kN+1}\binom{(k+1)N}{N}\\* &
	-\sum_{i=1}^{s-1}(s-i)\frac1{ki+1}\binom{(k+1)i}{i}\frac{k}{N-i+1}\binom{(k+1)(N-i)}{N-i}.	
\end{align*}
\begin{theorem}
		The generating funtion $\textsc{min}(z,w)$ where the coefficient of $z^{(k+1)N}w^s$ refers to the cumulative levels of the $s$-th \textsc{max}-turn, is given by
\begin{align*}
	\textsc{max}(z,w)
	&=	\frac{wk \uhat}{z(1-w)^2}
	-	\frac{wk\ubar }{z(1-w)^2}
	+\frac{w (\ubar-z)}{z^2\uhat^{k}(1-w)^2}
	-\frac{w^2 \ubar }{z(1-w)^2}.
\end{align*}
The coefficient  	$[w^sz^{(k+1)N}]	\textsc{max}(z,w)$ is for $N\ge s$ given by
\begin{align*}
\frac{sk}{kN+1}\binom{(k+1)N}{N}
	-\sum_{i=1}^{s-1}(s-i)\frac1{ki+1}\binom{(k+1)i}{i}\frac{k}{N-i+1}\binom{(k+1)(N-i)}{N-i}.	
\end{align*}
\end{theorem}

\section{The oscillation}

The cumulative function of the $s$-th oscillation (= total length of the $s$-th wavy line, cf.~Figure \ref{basic-fig} in all paths of length $(k+1)N$) is
\begin{equation*}
\textsc{osc}(z,w):=\textsc{max}(z,w)-\textsc{min}(z,w).
\end{equation*}

The simplication is easy:
\begin{align*}
\textsc{osc}(z,w)&=	\frac{wk \uhat}{z(1-w)^2}
-	\frac{wk\ubar }{z(1-w)^2}
+\frac{w (\ubar-z)}{z^2\uhat^{k}(1-w)^2}
-\frac{w^2 \ubar }{z(1-w)^2}\\
&-\bigg[\frac{kw\uhat}
{z(1-w)^2}+
\frac{(\ubar-z)}
{z^2(1-w)^2\uhat^{k}}-
\frac{(k+1)\overline{u}w}
{z(1-w)^2}\bigg]\\
&=\frac{w (\ubar-z)}{z^2\uhat^{k}(1-w)^2}
-\frac{(\ubar-z)}
{z^2\uhat^{k}(1-w)^2}-\frac{w^2 \ubar }{z(1-w)^2}+
\frac{w\overline{u}}
{z(1-w)^2}\\
&=		\frac{w\overline{u}}	{z(1-w)}-\frac{ (\ubar-z)}{z^2\uhat^{k}(1-w)}.
	\end{align*}
Further, by reading off coefficients ($s\ge1$),
\begin{align*}
	[w^s]&\textsc{osc}(z,w)=	
	\sum_{i=1}^s[w^i]\frac{\overline{u}}	{z}-\sum_{i=0}^s[w^i]\frac{ (\ubar-z)}{z^2\uhat^{k}}\\
	&=		\sum_{i=1}^sz^{(k+1)i}\frac1{kN+1}\binom{(k+1)i}{i}	-\frac1{z\uhat^{k}}\sum_{i=0}^s[w^i]\frac{ (\ubar-z)}{z}\\
		&=		\sum_{i=1}^sz^{(k+1)i}\frac1{kN+1}\binom{(k+1)i}{i}	-\frac1{z\uhat^{k}}\sum_{i=1}^s[w^i]\sum_{N\ge1}w^Nz^{(k+1)N}\frac1{kN+1}\binom{(k+1)N}{N}\\
		&=		\sum_{i=1}^sz^{(k+1)i}\frac1{kN+1}\binom{(k+1)i}{i}	-\frac1{z\uhat^{k}}\sum_{i=1}^sz^{(k+1)i}\frac1{ki+1}\binom{(k+1)i}{i}.
			\end{align*}
The next step is to read off $[w^sz^{(k+1)N}]\textsc{osc}(z,w)$ (for $N\ge s$):
\begin{align*}
	[w^s&z^{(k+1)N}]\textsc{osc}(z,w)
=		[z^{(k+1)N}]\sum_{i=1}^sz^{(k+1)i}\frac1{kN+1}\binom{(k+1)i}{i}\\&	\quad-[z^{(k+1)N}]\frac1{z\uhat^{k}}\sum_{i=1}^s[w^i]z^{(k+1)i}\frac1{ki+1}\binom{(k+1)i}{i}\\
	&=			-[z^{(k+1)N}]\bigg[z^{-k-1}-\sum_{\lambda\ge0}\frac{k}{\lambda+1}\binom{(k+1)\lambda}{\lambda}z^{(k+1)\lambda}\bigg]\sum_{i=1}^sz^{(k+1)i}\frac1{ki+1}\binom{(k+1)i}{i}\\
		&=	[z^{(k+1)N}]\sum_{i=1}^sz^{(k+1)i}\frac1{ki+1}\binom{(k+1)i}{i}\sum_{\lambda\ge0}\frac{k}{\lambda+1}\binom{(k+1)\lambda}{\lambda}
		z^{(k+1)\lambda}\\
		&=	\sum_{i=1}^sz^{(k+1)i}\frac1{ki+1}\binom{(k+1)i}{i}\frac{k}{N-i+1}\binom{(k+1)(N-i)}{N-i}.
\end{align*}

We summarize the results.
\begin{theorem}
	The cumulative (summed over all $k$-Dyck paths of length $(k+1)N$) difference of the $s$-th \textsc{max} turn and 
	the $(s+1)$-st \textsc{min} turn (nicknamed the $s$-th wavy line) is given by the explicit formula (for $N\ge s\ge1$)
\begin{align}\label{osc-formula}
	\sum_{i=1}^s\frac1{ki+1}\binom{(k+1)i}{i}\frac{k}{N-i+1}\binom{(k+1)(N-i)}{N-i}.
\end{align}
\end{theorem} 
In this paper his result that has been derived by manipulating power series only, was  obtained in \cite{AHS} by other (mostly bijective) methods.


	\section{Application: Hoppy walks}

Deng and Mansour \cite{deng}, compare also \cite{garden}, introduce a rabbit named Hoppy and let him move according to certain rules.
While the story about Hobby is charming and entertaining, we do not need this here and move straight ahead
to the enumeration issues. Eventually, the enumeration problem is one about $k$-Dyck paths ($k\ge1$).
The up-steps are $(1,k)$ and the down-steps are $(1,-1)$.

\begin{figure}[h]
	\begin{center}
		\begin{tikzpicture}[scale=0.5]
			\draw (0,0) -- (17,0);
			\draw (0,0) -- (0,10);
			\draw[thick](0,0)--(1,3);
			\draw[thick](5,6)--(6,9);
			\draw [decorate,decoration=snake,thick] (1,3) -- (5,6);
			\foreach \i in {0,...,8}
			{\draw[thick,red](6+\i,9-\i)--(7+\i,8-\i);
				\node[thick, red] at (6+\i,9-\i){$\bullet$};
			}
			\node[thick, red] at (6+5+4,9-5-4){$\bullet$};

			\node at (4,-1){$m$ up-steps};
			
			\node[thick] at (-1+0.4,9){$j$};

			\draw (-0.3,9) --(0.3,9);
		\end{tikzpicture}
	\end{center}
	\caption{The number of final down-steps}
\end{figure}
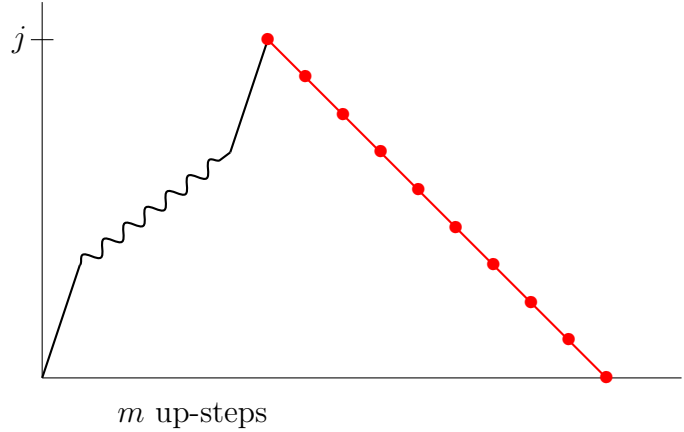	
The question is about the length of the sequence of down-steps printed in red. Or, phrased differently, how many $k$-Dyck paths end on level $j$, after $m$ up-steps, the last step being an up-step. The recent paper \cite{jcmcc} contains similar computations, although without the restriction that the last step must be an up-step. 

This question clearly belongs to the realm of the questions discussed here.

In fact, we are looking at the last oscillation (`wavy line'). If a $k$-Dyck path has length $(k+1)N$, the last oscillation occurs for $s=N$.
The formula (\ref{osc-formula}) leads for $s=N$ to
\begin{align}\label{osc-formula}
	\sum_{i=1}^N\frac1{ki+1}\binom{(k+1)i}{i}\frac{k}{N-i+1}\binom{(k+1)(N-i)}{N-i}.
\end{align}
Further considerations are in \cite{garden}.

	\section{Application: Sequences from the encyclopedia of integer sequences}

Now we move to instances of formula (\ref{osc-formula}) that appear in \cite{OEIS}. The principle reference is \cite{AHS}.
\begin{center}

\begin{tabular}{ c |c| l }
		$k=2$ & $s=N$ & A334680 \\
	$k=3$ & $s=N$ & A334682 \\
	$k=4$ & $s=N$ & A334719 \\
	\hline 
		$k=2$ & $s=1$ & A007226 \\
	$k=3$ & $s=1$ & A007228 \\
	$k=4$ & $s=1$ & A124724 \\
		\hline 
	$k=2$ & $s=2$ & A334640 \\
	$k=3$ & $s=2$ & A334645 \\
	$k=4$ & $s=2$ & A334646 \\
\end{tabular}
\end{center}

Tables of other similar sequences can be constructed easily with a computer algebra program, 
using explicit formul\ae\ in Theorems 1, 2, 3.

\clearpage

\bibliographystyle{plain}


\end{document}